\documentclass[11pt]{article}
\usepackage{theorem}
\usepackage{amssymb}
\usepackage{amsmath}
\usepackage{amsfonts}

\def \RR {\mathbb R}

\def \E {\mathcal E}
\def \F {\mathcal F}

\newtheorem{theorem}{Theorem}[section]
\newtheorem{lemma}[theorem]{Lemma}

\newtheorem{proposition}[theorem]{Proposition}
\newtheorem{corollary}[theorem]{Corollary}
 {\theorembodyfont{\rmfamily}}

\title{On John-Type Ellipsoids}
\author{B. Klartag\thanks{Supported in part by the Israel
Science Foundation, by the Minkowski Center for Geometry,
and by the European Research Training Network ``Analysis and Operators''.} \\
Tel Aviv University}
\date{}

\begin{document}

\maketitle

\abstract{Given an arbitrary convex symmetric body $K \subset
\RR^n$, we construct a natural and non-trivial continuous map
$u_K$ which associates ellipsoids to ellipsoids, such that the
L\"{o}wner-John ellipsoid of $K$ is its unique fixed point. A new
characterization of the L\"{o}wner-John ellipsoid is obtained, and we also
gain information regarding the contact points of inscribed ellipsoids with $K$.
}

\section{Introduction}
We work in $\RR^n$, yet we choose no canonical scalar product. A
centrally-symmetric ellipsoid in $\RR^n$ is any set of the form
$$ \left \{ \sum_{i=1}^n \lambda_i u_i \ ; \ \sum_i \lambda_i^2
 \leq 1, \ u_1,..,u_n \in \RR^n \right \}. $$
If $u_1,..,u_n$ are linearly independent, the ellipsoid is
non-degenerate. Whenever we mention an ``ellipsoid'' we mean a
centrally-symmetric non-degenerate one. Given an ellipsoid $\E
\subset \RR^n$, denote by $\langle \cdot , \cdot \rangle_{\E}$ the
unique scalar product such that $\E = \{ x \in \RR^n ; \langle x,
x \rangle_{\E} \leq 1 \}$. There is a group $O(\E)$ of linear
isometries of $\RR^n$ (with respect to the metric induced by
$\langle \cdot, \cdot \rangle_{\E}$), and a unique probability
measure $\mu_{\E}$ on $\partial \E$ which is invariant under
$O(\E)$. A body in $\RR^n$ is a centrally-symmetric convex set
with a non-empty interior. Given a body $K \subset \RR^n$, denote
by $\| \cdot \|_K$ the unique norm on $\RR^n$ such that $K$ is its
unit ball:
$$ \| x \|_K = \inf \{ \lambda > 0 ; x \in \lambda K \}. $$
Given a body $K \subset \RR^n$ and an ellipsoid $\E \subset
\RR^n$, denote
$$ M_{\E}^2(K) = M_{\E}^2(\| \cdot \|_K) = \int_{\partial \E} \| x \|_K^2 d\mu_{\E}(x). $$
This quantity is usually referred to as $M_2$. Let us consider the
following parameter:
\begin{equation}
J_K(\E) = \inf_{\F \subset K} M_{\E}(\F)
\label{minimizing_M_def}
\end{equation}
where the infimum runs over all ellipsoids $\F$ that are contained
in $K$. Since the set of all ellipsoids contained in $K$
(including degenerate ellipsoids) is a compact set with respect to
the Hausdorff metric, the infimum is actually attained. In
addition, the minimizing ellipsoids must be non-degenerate, since
otherwise $J_K(\E) = \infty$ which is impossible for a
body $K$. We are not so much interested in the
exact value of $J_K(\E)$, as in the ellipsoids where the minimum
is obtained.

\medskip In Section \ref{uniqueness} we prove that there exists a
unique ellipsoid for which
the minimum in (\ref{minimizing_M_def}) is attained. We shall
denote this unique ellipsoid by $u_K(\E)$, and we show that the
map $u_K$ is continuous. A finite measure $\nu$
on $\RR^n$ is called $\E$-isotropic if for any $\theta \in \RR^n$,
$$ \int \langle x, \theta \rangle_{\E}^2 d \nu(x) = L_{\nu}^2
\langle \theta, \theta \rangle_{\E} $$ where $L_{\nu}$ does not
depend on $\theta$. One of the important properties of the map
$u_K$ is summarized in the following proposition, to be proved in
Section \ref{charachterization}.

\begin{proposition}
Let $K \subset \RR^n$ be a body, and let $\E$, $\F \subset \RR^n$
be ellipsoids such that $\F \subset K$. Then $\F = u_K(\E)$ if and
only if there exists an $\E$-isotropic measure $\nu$ supported on
$\partial \F \cap \partial K$. \label{main_result}
\end{proposition}

In particular, given any Euclidean structure (i.e. scalar product)
in $\RR^n$, there is always a unique ellipsoid contained in $K$
with an isotropic measure supported on its contact points with
$K$. This unexpected fact
leads to a connection with the L\"{o}wner-John ellipsoid
of $K$, which is the (unique) ellipsoid of maximal volume
contained in $K$. By the characterization of the L\"{o}wner-John
ellipsoid due to John \cite{J} and Ball \cite{Ball} (see also
\cite{GM}), $u_K(\E) = \E$ if and only if the ellipsoid $\E$ is
the L\"{o}wner-John ellipsoid of $K$. Thus, we obtain the following:

\begin{corollary}
Let $K \subset \RR^n$ be a body, and let $\E \subset K$ be an
ellipsoid such that for any ellipsoid $\F \subset K$,
$$ M_{\E}(\F) \geq 1. $$
Then $\E$ is the L\"{o}wner-John ellipsoid of $K$.
\end{corollary}

As a byproduct of our methods, we also obtain an extremality
property of the mean width of the L\"{o}wner-John ellipsoid
(Corollary \ref{John_M_*}). In Section \ref{injective} we show
that the body $K$ is determined by the map $u_K$. Further evidence
for the naturalness of this map is demonstrated in Section
\ref{equivalent}, where we discuss optimization problems similar
to the optimization problem in (\ref{minimizing_M_def}), and
discover connections with the map $u_K$.

\section{Uniqueness}
\label{uniqueness}

Let $D$ be a minimizing ellipsoid in (\ref{minimizing_M_def}). We
will show that it is the only minimizing ellipsoid. We write $|x|
= \sqrt{\langle x, x \rangle_{D}}$.  An equivalent definition of
$J_K(\E)$ is the following:
\begin{equation}
J_K^2(\E) = \min \left \{ \int_{\partial \E} | T^{-1}(x)|^2
d\mu_{\E}(x) \ ; \ \| T:l_2^n \rightarrow X_K \| \leq 1 \right \}
\label{minimizing_M_def2}
\end{equation}
where $X_K = (\RR^n, \| \cdot \|_K)$ is the normed space whose
unit ball is $K$, and where $l_2^n = (\RR^n, | \cdot |)$. The
definitions are indeed equivalent; $T(D)$ is the ellipsoid from
definition (\ref{minimizing_M_def}), as clearly $\| x \|_{T(D)} =
| T^{-1}(x) |$. Since $D$ is a minimizing ellipsoid, $Id$ is a
minimizing operator in (\ref{minimizing_M_def2}). Note that in
(\ref{minimizing_M_def2}) it is enough to consider linear
transformations which are self adjoint and positive definite with
respect to $\langle \cdot, \cdot \rangle_D$. Assume on the
contrary that $T$ is another minimizer, where $T \neq Id$ is a
self adjoint positive definite operator. Let $\{e_1,..,e_n\}$ be
an orthogonal basis of eigenvectors of $T$, and let
$\lambda_1,..,\lambda_n > 0$ be the corresponding eigenvalues.
Consider the operator $S = \frac{Id + T}{2}$. Then $S$ satisfies
the norm condition in (\ref{minimizing_M_def2}), and by the strict
convexity of the function $x \mapsto \frac{1}{x^2}$ on $(0,
\infty)$,
$$ \int_{\partial \E} |S^{-1}(x)|^2 d\mu_{\E}(x) =
\int_{\partial \E} \sum_{i=1}^n \left( \frac{1}{\frac{1 +
\lambda_i}{2}} \right)^2 \langle x, e_i \rangle_D^2 d\mu_{\E}(x)
$$
$$ < \int_{\partial \E} \sum_{i=1}^n \frac{1 + \left( \frac{1}{\lambda_i}
\right)^2}{2}  \langle x, e_i \rangle_D^2 d\mu_{\E}(x) $$ $$ =
\frac{\int_{\partial \E} |x|^2 d \mu_{\E}(x) + \int_{\partial \E}
|T^{-1}(x)|^2 d\mu_{\E}(x)}{2} = J_K^2(\E) $$ since not all the
$\lambda_i$'s equal one, in contradiction to the minimizing
property of $Id$ and $T$. Thus the minimizer is unique, and we may
define a map $u_K$ which matches to any ellipsoid $\E$, the unique
ellipsoid $u_K(\E)$ such that $u_K(\E) \subset K$ and  $J_K(\E) =
M_{\E}(u_K(\E))$. It is easily verified that for any linear
operator $T$, and $t \neq 0$,
\begin{equation}
u_{TK}(T\E) = T u_K(\E), \label{u_K_linear}
\end{equation}
$$ u_K(t \E) = u_K(\E). $$
The second property means that the map $u_K$ is actually defined
over the ``projective space'' of ellipsoids. Moreover, the image
of $u_K$ is naturally a ``projective ellipsoid'' rather than an
ellipsoid: If $\E$ and $t \E$ both belong to the image of $u_K$,
then $t = \pm 1$. Nevertheless, we still formally define $u_K$ as
a map that matches an ellipsoid to an ellipsoid, and not as a map
defined over the ``projective space of ellipsoids''.

\smallskip Let us establish the continuity of the map $u_K$. One can
verify that $M_{\E}(\F)$ is a continuous function of $\E$ and $\F$
(using an explicit formula as in (\ref{scalar_integral}), for
example). Fix a body $K \subset \RR^n$, and denote by $X$ the
compact space of all (possibly degenerate) ellipsoids contained in
$K$. Then $M_{\E}(\F): X \times X \rightarrow [0, \infty]$ is
continuous. The map $u_K$ is defined only on a subset of $X$, the
set of non-degenerate ellipsoids. The continuity of $u_K$ follows
from the following standard lemma.

\begin{lemma}
Let $X$ be a compact metric space, and $f:X \times X \rightarrow
[0, \infty]$ a continuous function. Let $Y \subset X$, and assume
that for any $y \in Y$ there exists a unique $g(y) \in X$ such
that
$$ \min_{x \in X} f(x, y) = f(g(y), y). $$
Then $g:Y \rightarrow X$ is continuous.
\end{lemma}

\emph{Proof:} Assume that $y_n \rightarrow y$ in $Y$. The function
$\min_{x \in X} f(x, y)$ is continuous, and therefore
$$ \min_{x \in X} f(x, y_n) = f(g(y_n), y_n) \stackrel{n \rightarrow \infty}{\longrightarrow}
f(g(y),y) = \min_{x \in X} f(x, y). $$
 Since $X \times X$ is compact, $f$ is uniformly continuous and
\begin{eqnarray*}
\lefteqn{ \left| f(g(y_n), y) - f(g(y),y) \right| } \\
& \leq & \left| f(g(y_n), y) - f(g(y_n), y_n) \right| + \left|
f(g(y_n), y_n) - f(g(y),y) \right| \stackrel{n \rightarrow \infty}
\longrightarrow 0.
\end{eqnarray*}
 Therefore, for any convergent subsequence
$g(y_{n_k}) \rightarrow z$, we must have $f(z,y) = f(g(y), y)$ and
by uniqueness $z=g(y)$. Since $X$ is compact, necessarily $g(y_n)
\rightarrow g(y)$, and $g$ is continuous. \hfill $\square$

\section{Extremality conditions}
\label{charachterization}

There are several ways to prove the existence of the isotropic
measure announced in Proposition \ref{main_result}. One can adapt
the variational arguments from \cite{GM}, or use the Lagrange
multiplier technique due to John \cite{J} (as suggested by O.
Guedon). The argument we choose involves duality of linear
programming (see e.g. \cite{Bar}). For completeness, we state and
sketch the proof of the relevant theorem ($\langle \cdot, \cdot
\rangle$ is an arbitrary scalar product in $\RR^m$):

\begin{theorem}
Let $\{ u_{\alpha} \}_{\alpha \in \Omega} \subset \RR^m$, $\{
b_{\alpha} \}_{\alpha \in \Omega} \subset \RR$ and $c \in \RR^m$.
Assume that
$$ \langle x^0, c \rangle = \inf \{ \langle x, c \rangle ; \forall \alpha
\in \Omega, \langle x, u_{\alpha} \rangle \geq b_{\alpha} \} $$
and also $\langle x^0, u_{\alpha} \rangle \geq b_{\alpha}$ for any
$\alpha \in \Omega$. Then there exist $\lambda_1,..,\lambda_s > 0$
and $u_1,..,u_s \in \Omega^{\prime} = \{ \alpha \in \Omega ;
\langle x^0, u_{\alpha} \rangle = b_{\alpha} \} $ such that
$$ c = \sum_{i=1}^s \lambda_i u_i. $$
\label{lp_duality}
\end{theorem}
\vspace{-11pt}

\emph{Proof:} $K = \{ x \in \RR^m ; \forall \alpha \in \Omega,
\langle x, u_{\alpha} \rangle \geq b_{\alpha} \}$ is a convex
body. $x^0$ lies on its boundary, and $\{ x \in \RR^m ; \langle x,
c \rangle = \langle x^0, c \rangle \}$ is a supporting hyperplane
to $K$ at $x^0$. The vector $c$ is an inner normal vector to $K$
at $x^0$, hence $-c$ belongs to the cone of outer normal vectors
to $K$ at $x^0$. The crucial observation is that this cone is
generated by $-\Omega^{\prime}$ (e.g. Corollary 8.5 in chapter II
of \cite{Bar}), hence
$$ c \in \left \{ \sum_{i=1}^s \lambda_i u_i \ ; \ \forall i \  u_i \in
\Omega^{\prime}, \lambda_i \geq 0 \right \}. $$ \hfill $\square$

Let $K \subset \RR^n$ be a body, and let $\E \subset \RR^n$ be an
ellipsoid. This ellipsoid induces a scalar product in the space of
operators: if $T,S:\RR^n \rightarrow \RR^n$ are linear operators,
and $\{e_1,..,e_n\} \subset \RR^n$ is any orthogonal basis (with
respect to $\langle \cdot, \cdot \rangle_{\E}$), then
$$ \langle T, S \rangle_{\E} = \sum_{i,j} T_{i,j} S_{i,j}
$$
where $T_{i,j} = \langle Te_i, e_j \rangle_{\E}$ and $S_{i,j} =
\langle Se_i, e_j \rangle_{\E}$ are the entries of the
corresponding matrix representations of $T$ and $S$. This scalar
product does not depend on the choice of the orthogonal basis. If
$\F = \{ x \in \RR^n ; \langle x, Tx \rangle_{\E} \leq 1 \}$ is
another ellipsoid, then

\begin{eqnarray}
  \lefteqn{
M_{\E}^2(\F) = \int_{\partial \E} \langle x, Tx \rangle_{\E}
d\mu_{\E}(x)  \label{scalar_integral}}\\
& = & \sum_{i,j=1}^n T_{i,j} \int_{\partial \E} \langle x, e_i
\rangle_{\E} \langle x, e_j \rangle_{\E} d\mu_{\E}(x) =
\sum_{i,j=1}^n T_{i,j} \frac{\delta_{i,j}}{n} = \frac{1}{n}
\langle T, Id \rangle_{\E}.   \nonumber
\end{eqnarray}

The ellipsoid $\F = \{ x \in \RR^n ; \langle x, Tx \rangle_{\E}
\leq 1 \}$ is contained in $K$ if and only if for any $x \in
\partial K$,
$$ \langle x, Tx \rangle_{\E} = \langle x \otimes x, T \rangle_{\E} \geq 1 $$ where $(x \otimes
x)(y) = \langle x, y \rangle_{\E} x$ is a linear operator.
Therefore, the optimization problem (\ref{minimizing_M_def}) is
equivalent to the following problem:
$$ n J_K^2(\E) = \min \{ \langle T, Id \rangle_{\E} \ ; \ T \ is \ \E \! \scriptsize{-} \! positive, \ \forall x
\in \partial K \ \langle x \otimes x, T \rangle_{\E} \geq 1 \} $$
where we say that $T$ is $\E$-positive if it is self adjoint and
positive definite with respect to $\langle \cdot, \cdot
\rangle_{\E}$. Actually, the explicit positivity requirement is
unnecessary. If $K$ is non-degenerate and $\forall x \in \partial
K \ \langle T, x \otimes x \rangle_{\E} \geq 1$ then $T$ is
necessarily positive definite with respect to $\E$. This is a
linear optimization problem, in the space $\RR^m = \RR^{n^2}$. Let
$T$ be the unique self adjoint minimizer, and let $\F = \{ x \in
\RR^n ; \langle x, Tx \rangle_{\E} \leq 1 \}$ be the corresponding
ellipsoid. By Theorem \ref{lp_duality}, there exist
$\lambda_1,..,\lambda_s > 0$ and vectors $u_1,..,u_s \in \partial
K$ such that
\begin{enumerate}
\item For any $1 \leq i \leq s$ we have $\langle u_i \otimes u_i,
T \rangle_{\E} = 1$, i.e. $u_i \in \partial K \cap \partial \F$.
\item $Id = \sum_{i=1}^s \lambda_i u_i \otimes u_i$. Equivalently, for any $\theta \in \RR^n$,
$$ \sum_{i=1}^s \lambda_i \langle u_i, \theta \rangle_{\E}^2 = \langle \theta, \theta \rangle_{\E}^2. $$
\end{enumerate}

Hence we proved the following:
\begin{lemma} Let $K \subset \RR^n$ be a body and let $\E \subset \RR^n$ be an ellipsoid.
If $u_K(\E) = \F$, then there exist contact points $u_1,..,u_s \in
\partial K \cap
\partial \F$ and positive numbers $\lambda_1, .. , \lambda_s$ such that for any $\theta \in
\RR^n$,
$$
\sum_{i=1}^s \lambda_i \langle u_i, \theta \rangle^2_{\E} =
\langle \theta, \theta \rangle_{\E}.
$$
\label{decomposition}
\end{lemma}
\vspace{-11pt}

The following lemma completes the proof of Proposition
\ref{main_result}.
\begin{lemma}
Let $K \subset \RR^n$ be a body and let $\E, \F \subset \RR^n$ be
ellipsoids. Assume that $\F \subset K$ and that there exists a
measure $\nu$ supported on $\partial K \cap \partial \F$ such that
for any $\theta \in \RR^n$, $$ \int \langle x, \theta
\rangle_{\E}^2 d\nu(x) = \langle \theta, \theta \rangle_{\E}. $$
Then $u_K(\E) = \F$. \label{sufficient}
\end{lemma}

\emph{Proof:} Since $\int x \otimes x d\nu(x) = Id$, for any
operator $T$,
\begin{equation}
\int \langle Tx, x \rangle_{\E} d\nu(x) = \langle T, Id
\rangle_{\E} = n \int_{\partial \E} \langle Tx, x \rangle_{\E}
d\mu_{\E}(x). \label{equi_cond}
\end{equation}
where the last equality follows by (\ref{scalar_integral}). Let
$T$ be such that $\F = \{ x \in \RR^n ; \langle Tx, x \rangle_{\E}
\leq 1 \}$. By (\ref{equi_cond}),
$$ \nu(\partial \F) = \int_{\partial \F} \langle Tx, x
\rangle_{\E} d\nu(x) = n \int_{\partial \E} \langle Tx, x
\rangle_{\E} d\mu_{\E}(x) = n M_{\E}^2(\F). $$ Suppose that
$u_K(\E) \neq \F$. Then there exists a linear map $S \neq T$ such
that $\langle Sx, x \rangle_{\E} \geq 1$ for all $x \in
\partial K$ and such that
$$ \int_{\partial \E} \langle Sx, x \rangle_{\E} d\mu_{\E}(x) < M_{\E}^2(\F). $$
Therefore,
$$ \nu(\partial K) \leq \int_{\partial K} \langle Sx, x \rangle_{\E} d\nu(x)
 = n \int_{\partial \E} \langle Sx, x \rangle_{\E} d\mu_{\E}(x) < n M_{\E}^2(\F) $$
which is a contradiction, since $\nu(\partial K) = \nu(\partial \F)
= n M_{\E}^2(\F)$. \hfill $\square$

\medskip \emph{Remark:} This proof may be modified to provide an
alternative proof of John's theorem. Indeed, instead of minimizing
the linear functional $\langle T, Id \rangle$ we need to minimize
the concave functional $det^{1/n}(T)$. The minimizer still belongs
to the boundary, and minus of the gradient at this point belongs
to the normal cone.

\section{Different bodies have different maps}
\label{injective}

\begin{lemma}
Let $K, T \subset \RR^n$ be two closed bodies such that $T \not
\subset K$. Then there exists an ellipsoid $\F \subset T$ such
that $\F \not \subset K$ and $n$ linearly independent vectors
$v_1,..,v_n$ such that for any $1 \leq i \leq n$, $$ v_i \in
\partial \F \cap \partial C$$
where $C = conv(K, \F)$ and $conv$ denotes convex hull.
\label{ind_vecs}
\end{lemma}

\emph{Proof:} Let $U \subset T \setminus K$ be an open set whose
closure does not intersect $K$, and let $v_1^*,..,v_n^*$ be
linearly independent functionals on $\RR^n$ such that for any $y
\in K, \ z \in U$,
$$ v_i^*(y) < v_i^*(z) $$
for all $1 \leq i \leq n$. Let $\F \subset conv(U, -U)$ be an
ellipsoid that intersects $U$, and let $v_1,..,v_n \in \partial
\F$ be the unique vectors such that for $1 \leq i \leq n$,
$$ v_i^*(v_i) = \sup_{v \in \F} v_i^*(v). $$
Then $v_1,..,v_n$ are linearly independent. Also, $v_i^*(v_i) =
\sup_{v \in C} v_i^*(v)$ and hence $v_1,..,v_n$ belong to the
boundary of $C = conv(K, \F)$. \hfill $\square$

\begin{theorem}
Let $K, T \subset \RR^n$ be two closed bodies, such that for any
ellipsoid $\E \subset \RR^n$ we have $J_K(\E) = J_T(\E)$. Then $K = T$.
\label{J_equal}
\end{theorem}

\emph{Proof:} Assume $K \neq T$. Without loss of generality, $T
\not \subset K$. Let $\F$ and $v_1,..,v_n$ be the ellipsoid and
vectors from Lemma \ref{ind_vecs}. Consider the following bodies:
$$ L = conv \{\F, K \cap T \}, \ \ \ C = conv \{\F, K \}. $$
Then $v_1,..,v_n \in \partial C$ and also $v_1,..,v_n \in
\partial L$. Let $\langle \cdot, \cdot \rangle_{\E}$ be the scalar
product with respect to which these vectors constitute an
orthonormal basis. Then the uniform measure on $\{v_1,..,v_n\}$ is
$\E$-isotropic. By Proposition \ref{main_result}, $u_L(\E) =
u_C(\E) = \F$, and
\begin{equation}
J_L(\E) = M_{\E}(u_L(\E)) = M_{\E}(u_C(\E)) = J_C(\E).
\label{J_cond1}
\end{equation}
Since $K \subset C$, also $u_K(\E) \subset C$. Since $\F = u_C(\E)
\not \subset K$, we have $\F \neq u_K(\E)$. By the uniqueness of
the minimizing ellipsoid for $J_C(\E)$,
\begin{equation}
J_C(\E) = M_{\E}(\F) < M_{\E}(u_K(\E)) = J_K(\E). \label{C_K}
\end{equation}
 Since $L \subset T$ we have $J_T(\E) \leq J_L(\E)$. Combining
this with (\ref{J_cond1}) and (\ref{C_K}), we get
$$ J_T(\E) \leq J_L(\E) = J_C(\E) < J_K(\E) $$
and therefore $J_K(\E) \neq J_T(\E)$. \hfill $\square$

\begin{corollary}
Let $K, T \subset \RR^n$ be two closed bodies such that $u_K =
u_T$. Then $K = T$.
\end{corollary}

\emph{Proof:} For any ellipsoid $\E \subset \RR^n$,
$$ J_K(\E) = M_{\E}(u_K(\E)) = M_{\E}(u_T(\E)) = J_T(\E)
$$
and the corollary follows from Theorem \ref{J_equal}.

\section{Various optimization problems}
\label{equivalent}

Given an ellipsoid $\E \subset \RR^n$ and a body
$K \subset \RR^n$, define
$$ K_{\E}^{\circ} = \{ x \in \RR^n ;
\forall x \in K ,\ \langle x, y \rangle_{\E} \leq 1 \} $$ and also
$M^*_{\E}(K) = M_{\E} \left( K_{\E}^{\circ} \right)$. Consider the
following optimization problem:
\begin{equation}
\inf_{K \subset \F} M^*_{\E}(\F) \label{M_*_outside}
\end{equation}
where the infimum runs over all ellipsoids that contain $K$. Then
(\ref{M_*_outside}) is simply the dual, equivalent formulation of
problem (\ref{minimizing_M_def}) that was discussed above. Indeed, $\F$
is a minimizer in (\ref{M_*_outside}) if and only if $u_{K_{\E}^{\circ
}}(\E) = \F_{\E}^{\circ}$. An apriori different optimization
problem is the following:
\begin{equation}
I_K(\E) = \sup_{K \subset \F} M_{\E}(\F) \label{M_outside}
\end{equation}
where the supremum runs over all ellipsoids that contain $K$. The
characteristics of this problem are indeed different. For
instance, the supremum need not be attained, as shown by the
example of a narrow cylinder (in which there is a maximizing
sequence of ellipsoids that tends to an infinite cylinder) and
need not be unique, as shown by the example of a cube (any
ellipsoid whose axes are parallel to the edges of the cube, and
that touches the cube - is a maximizer. See also the proof of
Corollary \ref{John_M_*}). Nevertheless, we define
$$ \bar{u}_K(\E) = \left \{ \F \subset \RR^n \ ; \ \F \ is \ an \
ellipsoid, \ K \subset \F, \ \ M_{\E}(\F) = I_K(\E) \right \}.
$$
The dual, equivalent formulation of (\ref{M_outside}) means to
maximize $M^*$ among ellipsoids that are contained in $K$.
Apriori, $u_K(\E)$ and $\bar{u}_K(\E)$ do not seem to be related.
It is not clear why there should be a connection between
minimizing $M$ and maximizing $M^*$ among inscribed ellipsoids.
The following proposition reveals a close relation between the two
problems.

\begin{proposition}
Let $K \subset \RR^n$ be a body, and let $\E, \F \subset \RR^n$ be
ellipsoids. Then
$$ \F \in \bar{u}_K(\E) \ \ \ \Longleftrightarrow \ \ \ u_{K_{\F}^{\circ}}(\E) = \F. $$
 \label{main_result_outside}
\end{proposition}
\vspace{-11pt}

\emph{Proof:} \\

$\Longrightarrow$: We write $\F = \{ x \in \RR^n; \langle x, Tx
\rangle_{\E} \leq 1 \}$ for an $\E$-positive operator $T$. Since
$\F \in \bar{u}_K(\E)$, the operator $T$ is a maximizer of
$$ n I_K^2(\E) = \max \{ \langle S, Id \rangle_{\E} \ ; \ S \in L(n), \ \forall x \in
\partial K \ 0 \leq \langle S, x \otimes x \rangle_{\E} \leq 1 \}.
$$
where $L(n)$ is the space of linear operators acting on $\RR^n$.
Note that the requirement $\langle S, x \otimes x \rangle_{\E}
\geq 0$ for any $x \in \partial K$ ensures that $S$ is
$\E$-non-negative definite. This is a linear optimization problem.
Following the notation of Theorem \ref{lp_duality}, we rephrase
our problem as follows:
$$ -n I_K^2(\E) = \inf \{ \langle S, -Id \rangle_{\E} ; \forall x
 \in \partial K, \ \langle S, x \otimes x \rangle_{\E} \geq 0,
 \langle S, -x \otimes x \rangle_{\E} \geq -1 \}. $$
According to Theorem \ref{lp_duality}, since $T$ is a maximizer,
there exist $\lambda_1,..,\lambda_s
> 0$ and vectors $u_1,..,u_t,u_{t+1},..,u_s \in \partial K$ such that
\begin{enumerate}
\item For any $1 \leq i \leq t$ we have $\langle T, -u_i \otimes u_i,
 \rangle_{\E} = -1$, i.e. $u_i \in \partial K \cap \partial \F$.
\\
For any $t+1 \leq i \leq s$ we have $\langle T, u_i \otimes u_i
\rangle_{\E} = 0$.
 \item $Id = \sum_{i=1}^t \lambda_i u_i \otimes u_i -
 \sum_{i=t+1}^s \lambda_i u_i \otimes u_i$.
\end{enumerate}
Since we assumed that $\F$ is an ellipsoid, $T$ is $\E$-positive,
and it is impossible that $\langle T u_i, u_i \rangle_{\E} = 0$.
Hence, $t=s$ and there exists an $\E$-isotropic measure supported
on $\partial K \cap \partial \F$. Since $K \subset \F$, then $\F
\subset K_{\F}^{\circ}$ and $\partial K_{\F}^{\circ} \cap \partial
\F = \partial K \cap
\partial \F$. Therefore, there exists an $\E$-isotropic measure
supported on $\partial K_{\F}^{\circ} \cap \partial \F$. Since $\F
\subset K_{\F}^{\circ}$ we must have $u_{K_{\F}^{\circ}}(\E)=\F$,
according to Proposition \ref{main_result}.

\smallskip
$\Longleftarrow$: Since $u_{K_{\F}^{\circ}}(\E)=\F$, then $K \subset \F$ and we can write
$Id = \int x \otimes x d \nu(x)$ where $supp(\nu) \subset \partial K_{\F}^{\circ} \cap \partial \F = \partial K \cap \partial \F$. Reasoning as in Lemma \ref{sufficient}, $\nu(\partial K) = n M_{\E}^2(\F)$ and for any admissible operator $S$,
$$ \langle S, Id \rangle_{\E} = \int \langle x, Sx \rangle_{\E} d
\nu(x) \leq \nu(\partial K) = n M_{\E}^2(\F) $$ since $\langle x,
Sx \rangle_{\E} \leq 1$ for any $x \in \partial K$. Hence $T$ is a maximizer,
and $\F \in \bar{u}_K(\E)$. \hfill $\square$

\smallskip If $\E, \F \subset \RR^n$ are ellipsoids, then $K_{\E}^{\circ}$ is
a linear image of $K_{\F}^{\circ}$. By (\ref{u_K_linear}), the map
$u_{K_{\E}^{\circ}}$ is completely determined by $u_{K_{\F}^{\circ}}$.
Therefore, the family of maps $\{ u_{K_{\F}^{\circ}} ; \F \ is \ an
\ ellipsoid \}$ is determined by a single map $u_{K_{\E}^{\circ}}$,
for any ellipsoid $\E$. By Proposition \ref{main_result_outside},
this family of maps determines $\bar{u}_K$. Therefore, for any
ellipsoid $\E$, the map $u_{K_{\E}^{\circ}}$ completely determines
$\bar{u}_K$. Additional consequence of Proposition
\ref{main_result_outside} is the following:
\begin{corollary}
Let $K \subset \RR^n$ be a body, and
let $\E$ be its L\"{o}wner-John ellipsoid. Then for any ellipsoid
$\F \subset K$,
$$ M^*_{\E}(\F) \leq 1. $$
Equality occurs for $\F = \E$, yet there may be additional cases
of equality.
 \label{John_M_*}
\end{corollary}

\emph{Proof:} If $\E$ is the L\"{o}wner-John ellipsoid, then
$u_K(\E) = \E$. By Proposition \ref{main_result_outside}, $\E \in
\bar{u}_{K_{\E}^{\circ}}(\E)$. Dualizing, we get that $\E \subset
K$, and
$$ 1 = M_{\E}^*(\E) = \sup_{\F \subset K} M_{\E}^*(\F) $$
where the supremum is over all ellipsoids contained in $K$. This
proves the inequality. To obtain the remark about the equality
cases, consider the cross-polytope $K = \{ x \in \RR^n;
\sum_{i=1}^n |x_i| \leq 1 \}$, where $(x_1,..,x_n)$ are the
coordinates of $x$. By symmetry arguments, its L\"{o}wner-John
ellipsoid is $D = \{ x \in \RR^n; \sum_i x_i^2 \leq \frac{1}{n}
\}$. However, for any ellipsoid of the form
$$ \E = \left \{ x \in \RR^n ; \sum_i \frac{x_i^2}{\lambda_i} \leq 1 \right \}
$$ where the $\lambda_i$ are positive and $\sum_i \lambda_i = 1$,
we get that $\E \subset K$, yet $M_D^*(\E) = M_D^*(D) = 1$. \hfill
$\square$

\medskip \emph{Remarks:}
\begin{enumerate}
\item
 If $K \subset \RR^n$ is smooth and strictly convex,
then $\bar{u}_K$ is always a singleton. Indeed, if $K$ is strictly
convex and is contained in an infinite cylinder, it is also
contained in a subset of that cylinder which is an ellipsoid,
hence the supremum is attained. From the proof of Proposition
\ref{main_result_outside} it follows that if $\F_1,\F_2$ are
maximizers, then there exists an isotropic measure supported on
their common contact points with $K$. Since $K$ is smooth, it has
a unique supporting hyperplane at any of these contact points,
which is also a common supporting hyperplane of $\F_1$ and $\F_2$.
Since these common contact points span $\RR^n$, necessarily $\F_1
= \F_2$. Hence, if $K$ is smooth and strictly convex, only John
ellipsoid may cause an equality in Corollary \ref{John_M_*}. \item
If $\E$ is the L\"{o}wner-John ellipsoid of $K$, then for any
other ellipsoid $\F \subset K$ we have $M_{\E}(\F) > M_{\E}(\E) =
1$. This follows from our methods, yet it also follows immediately
from the fact that $\frac{1}{M_{\E}(\F)} \leq \left(
\frac{Vol(\F)}{Vol(\E)} \right)^{1/n}$, and from the uniqueness of
the L\"{o}wner-John ellipsoid.
\end{enumerate}

\bigskip \emph{Acknowledgement.} Part of the research was done
during my visit to the University of Paris $VI$, and I am very
grateful for their hospitality.

\end{document}